\documentclass[11pt]{amsart}
\usepackage{amssymb,amsmath}
\usepackage[all]{xy}

\title{Abstract homotopical methods for theoretical computer science}
\author[P. Gaucher]{Philippe Gaucher}
\address{Laboratoire PPS  (CNRS UMR 7126)\\ Universit{\'e} Paris 7--Denis Diderot\\
  Case 7014\\ 75205 PARIS Cedex 13 \\ France}
\email{gaucher@pps.jussieu.fr}
\urladdr{http://www.pps.jussieu.fr/{\~{}}gaucher/}

\newcommand{\de}{\partial}
\newcommand{\p}{\times}

\newcommand{\be}{\begin{equation}}
\newcommand{\ee}{\end{equation}}
\newcommand{\bea}{\begin{eqnarray}}
\newcommand{\eea}{\end{eqnarray}}
\newcommand{\beas}{\begin{eqnarray*}}
\newcommand{\eeas}{\end{eqnarray*}}

\newtheorem{thm}{Theorem}[section]
\newtheorem{prop}[thm]{Proposition}

\newtheorem{defn}[thm]{Definition}
\newtheorem{propdef}[thm]{Proposition and Definition}

\newcommand{\bd}{\begin{defn}}
\newcommand{\ed}{\end{defn}}
\newcommand{\bp}{\begin{prop}}
\newcommand{\ep}{\end{prop}}
\newcommand{\bth}{\begin{thm}}
\renewcommand{\eth}{\end{thm}}
\newcommand{\bpf}{\begin{proof}}
\newcommand{\epf}{\end{proof}}

\newcommand{\fl}[1]{\ar@{->}[ll]_-{#1}}
\newcommand{\fr}[1]{\ar@{->}[rr]^-{#1}}
\newcommand{\fd}[1]{\ar@{->}[dd]_-{#1}}
\newcommand{\fu}[1]{\ar@{->}[uu]^-{#1}}
\newcommand{\f}[2]{\ar@{->}[#1]|{#2}}
\newcommand{\ff}[2]{\ar@2{->}[#1]|{#2}}
\newcommand{\frr}[1]{\ar@{->}[rrrr]^-{#1}}

\newcommand{\ho}{{\mathbf{Ho}}}

\newcommand{\iso}{\cong}

\def\cartesien{%
  \ar@{-}[]+R+<6pt,-2pt>;[]+RD+<6pt,-6pt>%
  \ar@{-}[]+D+<2pt,-6pt>;[]+RD+<6pt,-6pt>%
}
\def\cocartesien{%
  \ar@{-}[]+L+<-6pt,+2pt>;[]+LU+<-6pt,+6pt>%
  \ar@{-}[]+U+<-2pt,+6pt>;[]+LU+<-6pt,+6pt>%
}
\def\hocartesien{%
  \ar@{-}[]+R+<6pt,-2pt>;[]+RD+<6pt,-6pt>_{h}%
  \ar@{-}[]+D+<2pt,-6pt>;[]+RD+<6pt,-6pt>%
}
\def\hococartesien{%
  \ar@{-}[]+L+<-6pt,+2pt>;[]+LU+<-6pt,+6pt>_{h}%
  \ar@{-}[]+U+<-2pt,+6pt>;[]+LU+<-6pt,+6pt>%
}

\newcommand{\brm}[1]{\rm{\mathbf{#1}}}

\DeclareMathOperator{\cons}{Diag}
\DeclareMathOperator{\hocons}{hoDiag}

\newcommand{\liminj}{\varinjlim}
\newcommand{\limproj}{\varprojlim}

\makeatletter

\def\varholim@#1#2{%
  \vtop{\m@th\ialign{##\cr
    \hfil$#1\operator@font holim$\hfil\cr
    \noalign{\nointerlineskip\kern1.5\ex@}#2\cr
    \noalign{\nointerlineskip\kern-\ex@}\cr}}%
}
\def\holimproj{%
  \mathop{\mathpalette\varholim@{\leftarrowfill@\textstyle}}\nmlimits@
}
\def\holiminj{%
  \mathop{\mathpalette\varholim@{\rightarrowfill@\textstyle}}\nmlimits@
}

\makeatother

\DeclareMathOperator{\id}{Id}

\DeclareMathOperator{\cell}{{\brm{cell}}}
\DeclareMathOperator{\cof}{{\brm{cof}}}
\DeclareMathOperator{\inj}{{\brm{inj}}}
\DeclareMathOperator{\Path}{{\rm{Path}}}
\DeclareMathOperator{\Cyl}{{\rm{Cyl}}}

\newcommand{\ddownarrow}{{\downarrow}}

\begin{document}

\begin{abstract} The purpose of this paper is to collect the
  homotopical methods used in the development of the theory of flows
  initialized by author's paper ``A model category for the homotopy
  theory of concurrency''. It is presented generalizations of the
  classical Whitehead theorem inverting weak homotopy equivalences
  between CW-complexes using weak factorization systems. It is also
  presented methods of calculation of homotopy limits and homotopy
  colimits using Quillen adjunctions and Reedy categories.
\end{abstract}

\maketitle

\tableofcontents

\section{Introduction}

The purpose of this paper is to collect the homotopical methods used
in the development of the theory of flows initialized by author's
paper \textit{A model category for the homotopy theory of concurrency}
\cite{model3} (see Table~\ref{recap}). An overview of the theory of
flows can be found in the two notes \cite{pgnote1} and \cite{pgnote2}.
The purpose of this paper is not to give a course in abstract homotopy
theory.  Indeed, there already exist several good introductions to
model category and more generally to abstract homotopy theory. For
model category, see the short papers \cite{MR1361887}
\cite{MR1916154}, or the books \cite{MR99h:55031} and
\cite{ref_model2}. For the relation between model category and
simplicial set, see the book \cite{MR2001d:55012}.  For cofibration
and fibration categories, see the book \cite{MR985099} or the notion
of ABC cofibration and fibration categories in \cite{ABChocolim}.  For
a more general setting allowing the development of the theory of
derived functors and homotopy limits and colimits without any model
category structure, see for example the books
\cite{monographie_hocolim} \cite{ref_model1}. The original reference
for model category is \cite{MR36:6480} but Quillen's axiomatization is
not used in this paper because it is obsolete. The Hovey's book
axiomatization is preferred.

It is also required a good knowledge of basic category theory.
Possible references are \cite{MR1712872} \cite{MR96g:18001a}
\cite{MR96g:18001b} and \cite{MR96g:18001c}.

The starting point is a complete cocomplete locally small category
$\mathcal{M}$ and a class of morphisms $\mathcal{W}$ (the \textit{weak
  equivalences}) modelling an equivalence relation between these
objects. One would like to consider them to be isomorphisms. It is
always possible to formally invert the weak equivalences by
considering the \textit{categorical localization}
$\mathcal{M}[\mathcal{W}^{-1}]$ of the category $\mathcal{M}$ with
respect to the morphisms of $\mathcal{W}$.  The categorical
localization is equipped with a canonical functor $\mathcal{M}
\rightarrow \mathcal{M}[\mathcal{W}^{-1}]$ which is the identity on
objects and which is universal for the property of taking all
morphisms of $\mathcal{W}$ to isomorphisms. The category
$\mathcal{M}[\mathcal{W}^{-1}]$ is difficult to understand since the
class of morphisms from $X$ to $Y$ is the quotient of the class of
finite zig-zag sequences $X=X_0 \longleftrightarrow X_1 \dots
\longleftrightarrow X_n=Y$, where the notation $A\longleftrightarrow
B$ means either a map from $A$ to $B$ or a map from $B$ to $A$, and
where all backward maps (i.e. pointing to the left) are weak
equivalences, divided by the equivalence relation generated by
removing or adding an identity map, and the identifications
$A\stackrel{f}\rightarrow B \stackrel{g}\rightarrow C =
A\stackrel{g\circ f}\longrightarrow C$, $A \stackrel{w}\rightarrow B
\stackrel{w}\leftarrow A= A \stackrel{\id_A}\longrightarrow A$ and $B
\stackrel{w}\leftarrow A \stackrel{w}\rightarrow B = B
\stackrel{\id_B}\longrightarrow B$ \cite{MR1712872} \cite{gz} or
\cite{ref_model1}. The class of morphisms in
$\mathcal{M}[\mathcal{W}^{-1}]$ between two objects needs even not be
a set. Two problems related to the study of
$\mathcal{M}[\mathcal{W}^{-1}]$ are treated in this paper.

First of all, it may be useful to construct notions of
\textit{cylinder} or \textit{cocylinder} functors related to
$\mathcal{W}$. This can be done by exhibiting a full subcategory
$\mathcal{M}_{good}$ of ``good'' objects of $\mathcal{M}$ such that
every object of $\mathcal{M}$ is weakly equivalent to a good object
and such that one has an isomorphism of categories
$\mathcal{M}_{good}/\!\!\sim \iso
\mathcal{M}_{good}[\mathcal{W}^{-1}]$ between the quotient of
$\mathcal{M}_{good}$ by a congruence $\sim$ generated by a cylinder or
a cocylinder functor and the categorical localization
$\mathcal{M}_{good}[\mathcal{W}^{-1}]$. Of course, this situation is
similar to the usual Whitehead statement inverting weak homotopy
equivalences between CW-complexes in classical algebraic topology
(\cite{MR1867354} Theorem~4.5). The main tool used for this problem is
the notion of \textit{weak factorization system}.

It may be also useful to calculate homotopy colimits and homotopy
limits with respect to $\mathcal{W}$. Let $\mathcal{B}$ be a small
category. Let $\mathcal{M}^\mathcal{B}$ be the category of functors
from $\mathcal{B}$ to $\mathcal{M}$. Let $\mathcal{W}_\mathcal{B}$ be
the class of morphisms $f:D \rightarrow E$ of
$\mathcal{M}^\mathcal{B}$ such that for every object $b$ of
$\mathcal{B}$, the map $f_b:D_b \rightarrow E_b$ is a weak
equivalence: let us call such a morphism an \textit{objectwise weak
  equivalence}. The constant diagram functor $\cons:\mathcal{M}
\rightarrow \mathcal{M}^\mathcal{B}$ is defined by $\cons(X)_b=X$ for
every object $b$ of $\mathcal{B}$ and $\cons(X)_\phi=\id_X$ (the
identity of $X$) for every morphism $\phi$ of $\mathcal{B}$. It has
both a left adjoint and a right adjoint since $\mathcal{M}$ is
complete and cocomplete, calculating respectively the colimit functor
and the limit functor from $\mathcal{M}^\mathcal{B}$ to $\mathcal{M}$.
It also induces a functor $\hocons:\mathcal{M}[\mathcal{W}^{-1}]
\rightarrow \mathcal{M}^\mathcal{B}[\mathcal{W}_\mathcal{B}^{-1}]$
between the localized categories because of the universal property
satisfied by $\mathcal{M}[\mathcal{W}^{-1}]$.  The problem is then to
give an explicit description of the left and right adjoints of
$\hocons$, if they exist. There is a considerable mathematical
literature about the subject. Connections between this problem and
\textit{model category} theory will be succinctly described. In
particular, the Reedy approach will be discussed. The last section is
devoted to the detailed description of several examples.

\begin{table}
\begin{center}
\begin{tabular}{|c|l|}
\hline
From this paper  & Used in ...\\
\hline
Theorem~\ref{W} & \cite{hocont} Theorem~3.27 and Theorem~4.6\\
\hline Proposition~\ref{rappel2} & \cite{exbranch} Theorem~5.5 ; \textit{\small see also  Theorem~\ref{Rcolim},\ref{Rlim}}\\
\hline Proposition~\ref{com} & \cite{exbranch} Lemma~8.6 and   Lemma~8.7\\
\hline Proposition~\ref{calcul} & \cite{3eme} Theorem~9.3, \cite{ccsprecub} Theorem~7.8\\
\hline Theorem~\ref{Reedymodel} & \cite{3eme} Theorem~8.4 ; \textit{\small see also Theorem~\ref{Rcolim},\ref{Rlim}}\\
\hline Theorem~\ref{Reedymodel} & \cite{4eme} Theorem~8.1, Theorem~8.2 and Theorem~8.3\\ 
\hline Theorem~\ref{Rcolim} & \cite{exbranch} Lemma~8.5, \cite{3eme} Theorem~7.5 \\
\hline Theorem~\ref{Rcolim} & \cite{4eme} Corollary~7.4, \cite{ccsprecub} Theorem~7.8\\
\hline Theorem~\ref{Rcolim} + left properness & \cite{3eme} Theorem~11.2, \cite{4eme} Theorem~9.1\\
\hline Theorem~\ref{Rlim} & \cite{model2} Theorem~IV.3.10 and  Theorem~IV.3.14\\
\hline Proposition~\ref{simp} & \textit{\small see Theorem~\ref{Rcolim}}\\
\hline
\end{tabular}
\end{center}
\caption{Summary of use of homotopical facts in homotopy theory of flows}
\label{recap}
\end{table}

\section{Whitehead theorem and weak factorization system}

Let $i:A\longrightarrow B$ and $p:X\longrightarrow Y$ be maps of the
category $\mathcal{M}$. Then $i$ has the {\rm left lifting property}
(LLP) with respect to $p$ (or $p$ has the {\rm right lifting property}
(RLP) with respect to $i$) if for every commutative square
\[
\xymatrix{
A\ar@{->}[dd]_{i} \ar@{->}[rr]^{\alpha} && X \ar@{->}[dd]^{p} \\
&&\\
B \ar@{-->}[rruu]^{g}\ar@{->}[rr]_{\beta} && Y,}
\]
there exists a morphism $g$ called a \textit{lift} making both
triangles commutative.

\bd \cite{MR2003h:18001} A {\rm weak factorization system} is a pair
$(\mathcal{L},\mathcal{R})$ of classes of morphisms of $\mathcal{M}$
such that the class $\mathcal{L}$ is the class of morphisms having the
LLP with respect to $\mathcal{R}$, such that the class $\mathcal{R}$
is the class of morphisms having the RLP with respect to $\mathcal{L}$
and such that every morphism of $\mathcal{M}$ factors as a composite
$r\circ \ell$ with $\ell\in \mathcal{L}$ and $r\in \mathcal{R}$. The
weak factorization system is {\rm functorial} if the factorization
$r\circ \ell$ can be made functorial.  \ed

If $K$ is a set of morphisms of $\mathcal{M}$, then the class of
morphisms of $\mathcal{M}$ that satisfy the RLP with respect to every
morphism of $K$ is denoted by $\inj(K)$ and the class of morphisms of
$\mathcal{M}$ that are transfinite compositions of pushouts of
elements of $K$ is denoted by $\cell(K)$. Denote by $\cof(K)$ the
class of morphisms of $\mathcal{M}$ that satisfy the LLP with respect
to every morphism of $\inj(K)$.

In a weak factorization system $(\mathcal{L},\mathcal{R})$, the class
$\mathcal{L}$ (resp.  $\mathcal{R}$) is completely determined by
$\mathcal{R}$ (resp.  $\mathcal{L}$). The class of morphisms
$\mathcal{L}$ is closed under composition, pushout, retract and binary
coproduct.  In particular, one has the inclusion $\cell(K)\subset
\cof(K)$. Dually, the class of morphisms $\mathcal{R}$ is closed under
composition, pullback, retract and binary product.

\bd Let $\mathcal{A}$ be a subcategory of $\mathcal{M}$. Then an
object $W$ is {\rm $\kappa$-small} relative to $\mathcal{A}$ for some
regular cardinal $\kappa$ if for every $\lambda$-sequence with
$\lambda\geq \kappa$
\[X_0 \rightarrow X_1\rightarrow X_2 \rightarrow \dots \rightarrow
X_\beta\rightarrow \dots \ \ (\beta<\lambda)\] such that $X_\beta
\rightarrow X_{\beta+1}$ is in $\mathcal{A}$ for every ordinal $\beta$
such that $\beta+1<\lambda$, the set map
$\liminj_{\beta<\lambda}\mathcal{M}(W,X_\beta) \rightarrow
\mathcal{M}(W,\liminj_{\beta<\lambda}X_\beta)$ is bijective. An object
$W$ is {\rm small} relative to $\mathcal{A}$ if it is $\kappa$-small
relative to $\mathcal{A}$ for some regular cardinal $\kappa$.  \ed

A set $K$ of morphisms of $\mathcal{M}$ \textit{permits the small
  object argument} if the domains of the morphisms of $K$ are small
relative to $\cell(K)$. In such a situation, the pair of classes of
morphisms $(\cof(K),\inj(K))$ can be made a functorial weak
factorization system using the small object argument
(\cite{ref_model2} Proposition~10.5.16, \cite{MR99h:55031}
Theorem~2.1.14). And moreover, every morphism of $\cof(K)$ is then a
retract of a morphism of $\cell(K)$ (\cite{MR99h:55031}
Corollary~2.1.15).

\bd \label{cg} A functorial weak factorization system
$(\mathcal{L},\mathcal{R})$ is {\rm cofibrantly generated} if there
exists a set $K$ of morphisms of $\mathcal{M}$ permitting the small
object argument such that $\mathcal{L}=\cof(K)$ and
$\mathcal{R}=\inj(K)$.  In this case, one explicitly supposes that the
functorial factorization is given by the small object argument
described in \cite{ref_model2} Proposition~10.5.16 or
\cite{MR99h:55031} Theorem~2.1.14, and not by any other method. \ed

Definition~\ref{cg} appears in \cite{MR1780498} in the context of
locally presentable category in the sense of \cite{MR95j:18001} as the
notion of \textit{small} weak factorization system. The reason of this
terminology is that in a locally presentable category, every
\textit{set} of morphisms permits the small object argument. Indeed,
every object of such a category is $\kappa$-presentable for a big
enough regular cardinal $\kappa$, and therefore $\kappa$-small
relative to the whole class of morphisms. So:

\bp  (\cite{MR1780498} Proposition~1.3)  
For every set of morphisms $K$ of a locally presentable category, the
pair of classes of morphisms $(\cof(K),\inj(K))$ is a cofibrantly
generated weak factorization system.
\ep

For the sequel, let us fix a functorial weak factorization system
$(\mathcal{L},\mathcal{R})$.

\bd \label{cyldef} Let $X$ be an object of $\mathcal{M}$. The {\rm
  cylinder object of $X$ with respect to $\mathcal{L}$} is the
functorial factorization
\[
\xymatrix{
X \oplus X \ar@{->}[rr]^-{\alpha(\id_X\oplus\id_X)}&& \Cyl_{\mathcal{L}}(X) \ar@{->}[rr]^-{\beta(\id_X\oplus\id_X)} && X}
\]
of the map $\id_X\oplus \id_X:X\oplus X \longrightarrow X$ by the
morphism $\alpha(\id_X\oplus\id_X):X\oplus X \longrightarrow
\Cyl_{\mathcal{L}}(X)$ of $\mathcal{L}$ composed with the morphism
$\beta(\id_X\oplus\id_X):\Cyl_{\mathcal{L}}(X) \longrightarrow X$ of
$\mathcal{R}$.  \ed

\bd An object $X$ of $\mathcal{M}$ is {\rm cofibrant with respect to
  $\mathcal{L}$} if the unique morphism $f_X:\varnothing
\longrightarrow X$, where $\varnothing$ is the initial object of
$\mathcal{M}$, is an element of $\mathcal{L}$.  Denote by
$\mathcal{M}_{cof}$ the category of cofibrant objects with respect to
$\mathcal{L}$.  \ed

\bd \label{leftdef} Let $f,g:X \rightrightarrows Y$ be two morphisms
of $\mathcal{M}$. A {\rm left homotopy with respect to} $\mathcal{L}$
from $f$ to $g$ is a morphism $H: \Cyl_{\mathcal{L}}X \longrightarrow
Y$ such that \[H\circ \alpha(\id_X\oplus\id_X) = f\oplus g.\] This
defines a reflexive and symmetric binary relation. The transitive
closure is denoted by $\sim^l_{\mathcal{L}}$.  \ed

It is clear that the weak factorization system
$(\mathcal{L},\mathcal{R})$ induces a weak factorization system
denoted in the same way on the full subcategory of cofibrant objects
with respect to $\mathcal{L}$. One then obtains the:

\bth \cite{ideeloc} \label{th1}
One has:
\begin{itemize}
\item The binary relation $\sim^l_{\mathcal{L}}$ does not depend on the
choice of the functorial factorization. 
\item The equivalence relation $\sim^l_\mathcal{L}$ is a congruence. 
\item Every object of $\mathcal{M}$ is isomorphic in $\mathcal{M}[\mathcal{R}^{-1}]$ to a 
cofibrant object with respect to
  $\mathcal{L}$.
\item The category $\mathcal{M}_{cof}/\!\!\sim^l_\mathcal{L}$ and
  $\mathcal{M}_{cof}[\mathcal{R}^{-1}]$ are isomorphic.
\item The category $\mathcal{M}_{cof}[\mathcal{R}^{-1}]$ is locally
  small.
\end{itemize}
\eth

The class of morphisms $\mathcal{R}$ plays the role of weak
equivalences in Theorem~\ref{th1}. This is exactly the situation
encountered by Y. Lafont and F. M{\'e}tayer in  \cite{MR1988395} 
\cite{lafont1} \cite{lafont2} \cite{metayer} in their study of 
higher dimensional rewriting systems using $\omega$-categories.

It is possible to dualize these results by working in the opposite
category $\mathcal{M}^{op}$ and by considering the weak factorization
system $(\mathcal{R}^{op},\mathcal{L}^{op})$. By definition, a
\textit{fibrant object of $\mathcal{M}$ with respect to $\mathcal{L}$}
is a cofibrant object of $\mathcal{M}^{op}$ with respect to
$\mathcal{R}^{op}$. The \textit{path object $\Path_{\mathcal{L}}X$} of
$X$ with respect to $\mathcal{L}$ is the cylinder object
$\Cyl_{\mathcal{L}}X$ of $X$ with respect to $\mathcal{R}^{op}$.
Finally, a \textit{right homotopy with respect to $\mathcal{L}$}
between two maps $f,g:X \rightrightarrows Y$ is a left homotopy
between $f^{op},g^{op}:Y \rightrightarrows X$ with respect to
$\mathcal{R}^{op}$. The transitive closure of this binary relation is
denoted by $\sim^r_{\mathcal{L}}$.  One obtains the following theorem,
in which the role of weak equivalences is now played by the morphisms
of $\mathcal{L}$:

\bth \cite{ideeloc}
One has:
\begin{itemize}
\item The binary relation $\sim^r_{\mathcal{L}}$ does not depend on the
choice of the functorial factorization. 
\item The equivalence relation $\sim^r_\mathcal{L}$ is a congruence. 
\item Every object of $\mathcal{M}$ is isomorphic in
  $\mathcal{M}[\mathcal{L}^{-1}]$ with a fibrant object with respect
  to $\mathcal{L}$.
\item The category $\mathcal{M}_{fib}/\!\!\sim^l_\mathcal{L}$ and
  $\mathcal{M}_{fib}[\mathcal{L}^{-1}]$ are isomorphic, where $
  \mathcal{M}_{fib}$ is the full subcategory of fibrant objects with
  respect to $\mathcal{L}$.
\item The category $\mathcal{M}_{fib}[\mathcal{L}^{-1}]$ is locally
  small.
\end{itemize}
\eth

The interest of the dual theorem is that it can be improved as
follows:

\bth\cite{hocont}\label{W} Suppose that the weak factorization system
$(\mathcal{L},\mathcal{R})$ is cofibrantly generated and that every
map of $\mathcal{L}$ is a monomorphism. Then the inclusion functor
$\mathcal{M}_{fib} \subset \mathcal{M}$ induces an equivalence of
categories $\mathcal{M}_{fib}/\!\!\sim^r_\mathcal{L} \simeq
\mathcal{M}[\mathcal{L}^{-1}]$. In particular, the category
$\mathcal{M}[\mathcal{L}^{-1}]$ is locally small.  \eth

Theorem~\ref{W} is proved using a fibrant replacement functor
$R_\mathcal{L}:\mathcal{M} \rightarrow \mathcal{M}_{fib}$ defined by
factoring the natural map $X \rightarrow \mathbf{1}$ from an object
$X$ to the terminal object as a composite $X \rightarrow
R_\mathcal{L}(X) \rightarrow \mathbf{1}$ using the weak factorization
system $(\mathcal{L},\mathcal{R})$. The factorization must be obtained
using the small object argument. It is then possible to prove that the
functor $R_\mathcal{L}$ is inverse to the inclusion functor up to
isomorphism of functors.

Theorem~\ref{W} is actually used in \cite{hocont} with $\mathcal{M}$
replaced by the full subcategory of cofibrant objects of a model
category in the sense of Definition~\ref{modelc}. It enables us to
prove a Whitehead theorem for the full dihomotopy relation on flows
(see \cite{1eme} for an informal introduction about dihomotopy of
flows). By \cite{nonexistence}, the dihomotopy relation on the
category of flows does not correspond to any model category structure
in the sense of Definition~\ref{modelc} but it can be fully described
using the weak S-homotopy model structure constructed in \cite{model3}
and an additional weak factorization system modelling refinement of
observation.

\section{Model category and Quillen adjunction}

It is introduced in this section the fundamental tool of \textit{model
  category} and \textit{Quillen derived functor}. The difference with
the preceding section is that we now have two weak factorization
systems interacting with each other and which are related to the class
of weak equivalences.

\bd \label{modelc}
A {\rm model category} is a complete and
cocomplete category $\mathcal{M}$ equipped with three classes of morphisms
$({\rm Cof},{\rm Fib},\mathcal{W})$ (resp. called the classes of
{\rm cofibrations}, {\rm fibrations} and {\rm weak equivalences}) such that:
\begin{enumerate}
\item the class of morphisms $\mathcal{W}$ is closed under retracts
  and satisfies the 2-out-of-3 property i.e.: if $f$ and $g$ are
  morphisms of $\mathcal{M}$ such that $g\circ f$ is defined and two
  of $f$, $g$ and $g\circ f$ are weak equivalences, then so is the
  third.
\item the pairs $({\rm Cof}\cap \mathcal{W}, {\rm Fib})$ and
$({\rm Cof}, {\rm Fib}\cap \mathcal{W})$ are both functorial weak
factorization systems.
\end{enumerate}
The triple $({\rm Cof},{\rm Fib},\mathcal{W})$ is called a {\rm model
  structure}. An element of ${\rm Cof}\cap \mathcal{W}$ is called a
{\rm trivial cofibration}. An element of ${\rm Fib}\cap \mathcal{W}$
is called a {\rm trivial fibration}. The categorical localization
$\ho(\mathcal{M}) := \mathcal{M}[\mathcal{W}^{-1}]$ is called the {\rm
  homotopy category} of $\mathcal{M}$.  The model category
$\mathcal{M}$ is {\rm cofibrantly generated} if both weak
factorization systems $({\rm Cof}\cap \mathcal{W}, {\rm Fib})$ and
$({\rm Cof}, {\rm Fib}\cap \mathcal{W})$ are cofibrantly generated.
\ed

It is introduced in \cite{model3} a cofibrantly generated model
structure for the study of concurrency.  The objects are called
\textit{flows}, they are actually categories without identity maps
enriched over compactly generated topological spaces (more details for
this kind of topological spaces in \cite{MR90k:54001}
\cite{MR2000h:55002}, the appendix of \cite{Ref_wH} and also the
preliminaries of \cite{model3}). The weak equivalences are the
morphisms of flows inducing a bijection between the sets of objects
and a weak homotopy equivalence between the spaces of morphisms; the
fibrations are the morphisms of flows inducing a fibration on the
space of morphisms; finally the cofibrations are determined by the
left lifting property with respect to trivial fibrations. Another
cofibrantly generated model category relevant for concurrency theory
is the model structure constructed by K. Worytkiewicz on cubical sets
\cite{W}. The common feature of the two model structures is that the
directed segment is not equivalent to a point.  It is very important
for the preservation of the causal structure of the underlying time
flow. See the long introduction of \cite{hocont} for further
explanations.

An object $X$ of a model category $\mathcal{M}$ is \textit{cofibrant}
(resp. \textit{fibrant}) if and only if it is cofibrant with respect
to ${\rm Cof}$ (resp. fibrant with respect to ${\rm Cof}\cap
\mathcal{W}$). The canonical morphism $\varnothing\rightarrow X$
functorially factors as a composite $\varnothing\rightarrow X^{cof}
\rightarrow X$ of a cofibration $\varnothing\rightarrow X^{cof}$
followed by a trivial fibration $X^{cof} \rightarrow X$.
Symmetrically, the canonical morphism $X\rightarrow \mathbf{1}$
functorially factors as a composite $X\rightarrow X^{fib} \rightarrow
\mathbf{1}$ of a trivial cofibration $X\rightarrow X^{fib}$ followed
by a fibration $X^{fib} \rightarrow \mathbf{1}$.

\bd The functor $X \mapsto X^{cof}$ is called the {\rm cofibrant
  replacement functor}.  The functor $X \mapsto X^{fib}$ is called the
{\rm fibrant replacement functor}.  \ed

\bp  (\cite{MR99h:55031} Theorem~1.2.10 and \cite{ref_model2} Theorem~8.3.9) 
Let $\mathcal{M}$ be a model category. Let $\mathcal{M}_{cof,fib}$ be
the full subcategory of cofibrant-fibrant objects of
$\mathcal{M}$. Then the inclusion functor $\mathcal{M}_{cof,fib}
\subset \mathcal{M}$ induces an equivalence of categories
$\mathcal{M}_{cof,fib}/\!\!\sim \simeq \ho(\mathcal{M})$ where the
congruence $\sim$ is left homotopy with respect to ${\rm Cof}$, or
equivalently right homotopy with respect to ${\rm Cof}\cap
\mathcal{W}$ (the two congruences coincide on cofibrant-fibrant
objects). In particular, the homotopy category $\ho(\mathcal{M})$ is
locally small. \ep

The model category $\mathcal{M}$ is \textit{left proper} (resp.
\textit{right proper}) if every pushout (resp. pullback) of a weak
equivalence along a cofibration (resp. fibration) is a weak
equivalence.  A model category $\mathcal{M}$ is \textit{proper} if it
is both left and right proper.

The following proposition explains the relation between properness,
fibrancy and cofibrancy. It can also be viewed as a method of
construction of weak equivalences.

\bp \label{rel} (Reedy) (\cite{ref_model2} Proposition~13.1.2) 
Let $\mathcal{M}$ be a model category. Then 
\begin{itemize}
\item Every pushout of a weak equivalence between cofibrant objects
  along a cofibration is a weak equivalence.
\item Every pullback of a weak equivalence between fibrant objects
  along a fibration is a weak equivalence.
\end{itemize}
\ep

The two consequences of Proposition~\ref{rel} are:
\begin{itemize}
\item A model category where all objects are cofibrant (like the model
  category of simplicial sets \cite{MR2001d:55012}) is left proper.
\item A model category where all objects are fibrant (like the model
  category of compactly generated topological spaces
  \cite{MR99h:55031} or the model category of flows \cite{model3}) is
  right proper.
\end{itemize}

The model categories of simplicial sets, of compactly generated
topological spaces and of flows are actually all of them proper, i.e.
both left and right proper.

\begin{propdef} \label{rappel2} (\cite{ref_model2} Proposition~8.5.3
  and Proposition~8.5.4) A Quillen adjunction is a pair of adjoint
  functors $F:\mathcal{M}\rightleftarrows \mathcal{N}:G$ between the
  model categories $\mathcal{M}$ and $\mathcal{N}$ such that one of
  the following equivalent properties holds:
\begin{enumerate}
\item $F$ preserves both cofibrations and trivial cofibrations.
\item $G$ preserves both fibrations and trivial fibrations.
\item $F$ preserves both cofibrations between cofibrant objects and
  trivial cofibrations (D. Dugger).
\item $G$ preserves both fibrations between fibrant objects and
  trivial fibrations (D. Dugger).
\end{enumerate}
One says that $F$ is a {\rm left Quillen functor}.  One says that $G$
is a {\rm right Quillen functor}. Moreover, any left Quillen functor
preserves weak equivalences between cofibrant objects and any right
Quillen functor preserves weak equivalences between fibrant objects.
\end{propdef}

The branching space functor from the category of flows to the category
of compactly generated topological spaces, defined in \cite{exbranch}
and studied in \cite{3eme} is an example of left Quillen functor if
the category of flows is equipped with the weak S-homotopy model
structure constructed in \cite{model3} and if the category of
compactly generated topological spaces is equipped with the usual
Quillen model structure.

Let $F:\mathcal{M}\rightleftarrows \mathcal{N}:G$ be a Quillen
adjunction between the model categories $\mathcal{M}$ and
$\mathcal{N}$. The cofibrant replacement functor takes weak
equivalence to weak equivalence by the 2-out-of-3 property. So the
functor $F\circ (-)^{cof}:\mathcal{M} \rightarrow \mathcal{N}$ induces
a unique functor $\mathbf{L}F:\ho(\mathcal{M}) \rightarrow
\ho(\mathcal{N})$ making the following diagram commutative
\[
\xymatrix{
\mathcal{M} \fr{F\circ (-)^{cof}}\fd{i_\mathcal{M}} && \mathcal{N} \fd{i_\mathcal{N}}\\
&& \\
\ho(\mathcal{M}) \fr{\mathbf{L}F} && \ho(\mathcal{N}).}
\] 
because of the universal property satisfied by $\ho(\mathcal{M})$. In
the same way, the fibrant replacement functor takes weak equivalence
to weak equivalence by the 2-out-of-3 property. So the functor $G\circ
(-)^{fib}:\mathcal{N} \rightarrow \mathcal{M}$ induces a unique
functor $\mathbf{R}G:\ho(\mathcal{N}) \rightarrow \ho(\mathcal{M})$
making the following diagram commutative
\[
\xymatrix{
\mathcal{N} \fr{G\circ (-)^{fib}}\fd{i_\mathcal{N}} && \mathcal{M} \fd{i_\mathcal{M}}\\
&& \\
\ho(\mathcal{N}) \fr{\mathbf{R}G} && \ho(\mathcal{M})}
\] 
because of the universal property satisfied by $\ho(\mathcal{N})$.

\bd The functor $\mathbf{L}F:\ho(\mathcal{M}) \rightarrow \ho(\mathcal{N})$  is called the 
{\rm total left derived functor} of $F$. The functor
$\mathbf{R}G:\ho(\mathcal{N}) \rightarrow \ho(\mathcal{M})$ is called the {\rm total
right derived functor} of $G$. \ed

The following theorem is the main tool for the sequel:

\bth \label{adj} (Quillen) (\cite{MR99h:55031} Lemma~1.3.10 or \cite{ref_model2}
Theorem~8.5.18) The pair of functors $\mathbf{L}F:\ho(\mathcal{M})
\rightleftarrows \ho(\mathcal{N}):\mathbf{R}G$ is a categorical
adjunction. \eth

The natural transformation $(-)^{cof} \Rightarrow \id$ induces a
natural transformation $\mathbf{L}F\circ i_\mathcal{M}=
i_\mathcal{N}\circ F \circ (-)^{cof} \Rightarrow i_\mathcal{N}\circ
F$. Of course, if $F$ already preserves weak equivalences, then one
has the isomorphism of functors $\mathbf{L}F\circ i_\mathcal{M}
\iso i_\mathcal{N}\circ F$. And the natural transformation $\id
\Rightarrow (-)^{fib}$ induces a natural transformation $i_\mathcal{M}
\circ G \Rightarrow i_\mathcal{M} \circ G \circ (-)^{fib} =
\mathbf{R}G \circ i_\mathcal{N}$. And if $G$ already preserves weak
equivalences, then one has the isomorphism of functors $i_\mathcal{M}
\circ G \iso \mathbf{R}G \circ i_\mathcal{N}$.

In fact, the functor $\mathbf{L}F:\ho(\mathcal{M}) \rightarrow
\ho(\mathcal{N})$ is the right Kan extension of the functor
$i_\mathcal{N}\circ F:\mathcal{M} \rightarrow \ho(\mathcal{N})$ along
the canonical functor $i_\mathcal{M}:\mathcal{M} \rightarrow
\ho(\mathcal{M})$ and the functor $\mathbf{R}G:\ho(\mathcal{N})
\rightarrow \ho(\mathcal{M}) $ is the left Kan extension of the
functor $i_\mathcal{M} \circ G:\mathcal{N} \rightarrow
\ho(\mathcal{M})$ along the canonical functor
$i_\mathcal{N}:\mathcal{N} \rightarrow \ho(\mathcal{N})$. So the
functor $\mathbf{L}F\circ i_\mathcal{M}$ is the closest approximation
of $F$ preserving weak equivalences and the functor $\mathbf{R}G \circ
i_\mathcal{N}$ is the closest approximation of $G$ preserving weak
equivalences.

If $F_1$ and $F_2$ are two composable left Quillen adjoints, then
$F_1\circ F_2$ is a left Quillen adjoint and the natural
transformation $F_1\circ (-)^{cof}\circ F_2 \circ (-)^{cof}
\Rightarrow F_1\circ F_2 \circ (-)^{cof}$ induces an isomorphism of functors 
$\mathbf{L}(F_1) \circ \mathbf{L}(F_2) \iso \mathbf{L}(F_1\circ F_2)$.
Similarly, if $G_1$ and $G_2$ are two composable right Quillen
adjoints, then $G_1\circ G_2$ is a right Quillen adjoint and the
natural transformation $G_1\circ G_2 \circ (-)^{fib} \Rightarrow
G_1\circ (-)^{fib}\circ G_2 \circ (-)^{fib}$ induces an isomorphism of
functors $\mathbf{R}(G_1\circ G_2)\iso \mathbf{R}(G_1) \circ
\mathbf{R}(G_2)$. See \cite{MR99h:55031} Theorem~1.3.7 for further details.

\section{Homotopy limit and homotopy colimit}

We want to give in this section constructions of homotopy limits and
homotopy colimits in the following situations:
\begin{itemize}
\item a construction of $\holiminj$ for every cofibrantly generated
  model category $\mathcal{M}$
\item a construction of $\holimproj$ for every cofibrantly generated
  model category $\mathcal{M}$ such that $\mathcal{M}$ is locally
  presentable~\footnote{ The construction requires that the class of
  weak equivalences satisfies the solution set
  condition. J. Rosick{\'y} has a proof that every locally presentable
  cofibrantly generated model category has an accessible class of weak
  equivalences. Thus in particular, it satisfies the solution set
  condition: these two conditions are not equivalent; a large cardinal
  axiom is needed for the converse \cite{MR1319969}. If $\mathcal{M}$
  is left proper and simplicial, then \cite{stablemodele}
  Proposition~3.2 provides an accessible functor from $\mathcal{M}$ to
  simplicial sets detecting weak equivalences. So by
  \cite{MR95j:18001} Corollary~2.45, the class of weak equivalences
  satisfies the solution set condition since the class of weak
  equivalences of simplicial sets is accessible (\cite{MR1780498}
  Example~3.1 and \cite{MR0491680} \cite{MR0491681}).}
\item a construction of $\holiminj$ when $\mathcal{B}$ is a Reedy
  category with fibrant constants
\item a construction of $\holimproj$ when $\mathcal{B}$ is a Reedy
  category with cofibrant constants. 
\end{itemize}
It is not always possible to get such a situation. Hence the interest
of other approaches like \cite{monographie_hocolim} \cite{ref_model1},
the simplicial technique of \cite{ref_model2} Chapter~18 and the
technique of frames of \cite{ref_model2} Chapter~19.

Before going further, it may be useful to point out that the homotopy
category of any model category is weakly complete and weakly
cocomplete.  Weak limit and weak colimit satisfy the same property as
limit and colimit except the unicity.  Weak small (co)products
coincide with small (co)products.  Weak (co)limits are constructed
using small (co)products and weak (co)equalizers in the same way as
(co)limits are constructed by small (co)products and (co)equalizers
(\cite{MR1712872} Theorem~1 p109). And a weak coequalizer 
\[A\stackrel{f,g}\rightrightarrows B \stackrel{h}\longrightarrow D\] 
is given by a weak pushout 
\[
\xymatrix{
B \fr{h} && D \\
&& \\
A\sqcup B \fu{(f,\id_B)}\fr{(g,\id_B)} && B. \fu{h}}
\] 
And finally, weak pushouts (resp. weak pullbacks) are given by
homotopy pushouts (resp. homotopy pullbacks) by
\cite{rep} Remark~4.1. See also \cite{MR920963} Chapter~III.

Let us come back now to homotopy limits and colimits. The principle,
exposed in this section, for calculating $\holiminj$ is the
construction of a model structure on the category of diagrams
$\mathcal{M}^\mathcal{B}$ such that the colimit functor
$\liminj:\mathcal{M}^\mathcal{B}
\rightarrow \mathcal{M}$ becomes a left Quillen functor.
Symmetrically, the principle for calculating $\holimproj$ is the
construction of a model structure on the category of diagrams
$\mathcal{M}^\mathcal{B}$ such that the limit functor
$\limproj:\mathcal{M}^\mathcal{B} \rightarrow \mathcal{M}$ becomes a
right Quillen functor. Indeed, if the categorical adjunction
$\liminj:\mathcal{M}^\mathcal{B} \rightleftarrows \mathcal{M}:\cons$
is a Quillen adjunction, then the natural transformation
$\hocons \Rightarrow \hocons\circ
\ho((-)^{fib})=\mathbf{R}\cons$ is an isomorphism of functors. So the
left adjoint of $\hocons$ exists by Theorem~\ref{adj} and is
the homotopy colimit.  And if the categorical adjunction
$\cons:\mathcal{M} \rightleftarrows\mathcal{M}^\mathcal{B} :\limproj$
is a Quillen adjunction, then the natural transformation
$\mathbf{L}\cons=\hocons\circ \ho((-)^{cof}) \Rightarrow
\hocons$ is an isomorphism of functors as well. So the right
adjoint of $\hocons$ exists by Theorem~\ref{adj} and is the homotopy
limit. A straightforward consequence is:

\bp \label{com} A total left derived functor commutes with homotopy
colimits.  A total right derived functor commutes with homotopy
limits. \ep

Proposition~\ref{com} is used in the proofs of \cite{exbranch}
Lemma~8.6 and \cite{exbranch} Lemma~8.7. The following proposition
gives an example of calculation of homotopy colimit, used in
\cite{3eme} Theorem~9.3 and in \cite{ccsprecub} Theorem~7.8:

\bp \label{calcul} (\cite{ref_model2} Proposition~18.1.6 and
Proposition~14.3.13) The homotopy colimit of a diagram of contractible
topological spaces over $\mathcal{B}$ is homotopy equivalent to the
classifying space of $\mathcal{B}$ (\cite{ref_model2} Chapter~14 or
\cite{MR0232393} \cite{MR0338129} for a definition of the classifying
space). In particular, if $\mathcal{B}$ has an initial or a terminal
object, then this homotopy colimit is contractible as well.  \ep

There exist two general theorems providing model structures on
$\mathcal{M}^\mathcal{B}$ such that the colimit functor (resp. the
limit functor) is a left (resp. right) Quillen functor. The following
theorem ensures the existence of homotopy colimit for any cofibrantly
generated model category $\mathcal{M}$ and for any small category
$\mathcal{B}$:

\bth \label{hocolim} (\cite{ref_model2} Theorem~11.6.1, Theorem~11.6.8
and Theorem~13.1.14) Let us suppose $\mathcal{M}$ equipped with a
cofibrantly generated model structure.  Then there exists a unique
model structure on $\mathcal{M}^\mathcal{B}$ such that the fibrations
are the objectwise fibrations and the weak equivalences the objectwise
weak equivalences.  Moreover, one has:
\begin{itemize}
\item Every cofibration of this model structure is an objectwise
cofibration. 
\item This model structure is cofibrantly generated. 
\item The colimit functor $\liminj:\mathcal{M}^\mathcal{B}
\rightarrow \mathcal{M}$ is a left Quillen functor.  
\item If $\mathcal{M}$ is left proper (resp. right proper, proper), then 
so is $\mathcal{M}^\mathcal{B}$. 
\end{itemize}
\eth

The following theorem ensures the existence of homotopy limit for any
cofibrantly generated model category $\mathcal{M}$ with $\mathcal{M}$
locally presentable with a class of weak equivalences satisfying the
solution set condition and for any small category $\mathcal{B}$:

\bth (Unknown reference) \label{holim} Let us suppose $\mathcal{M}$
locally presentable and cofibrantly generated.  Then there exists a
unique model structure on $\mathcal{M}^\mathcal{B}$ such that the
cofibrations are the objectwise cofibrations and the weak equivalences
the objectwise weak equivalences. Moreover, one has:
\begin{itemize}
\item Every fibration of this model structure is an objectwise
  fibration.
\item This model structure is cofibrantly generated.
\item The limit functor
$\limproj:\mathcal{M}^\mathcal{B} \rightarrow \mathcal{M}$ is a right
Quillen functor.
\end{itemize}
\eth

Theorem~\ref{holim} is a direct consequence of a theorem due to J.
Smith and exposed in \cite{MR1780498} Theorem~1.7.
Theorem~\ref{holim} is very close to the statement of \cite{hitopos}
Proposition~A.3.3.3.  Theorem~\ref{hocolim} and Theorem~\ref{holim}
can be applied to the case of category of simplicial presheaves
\cite{MR1300636}. The model structure of Theorem~\ref{hocolim} is then
known as the Bousfield-Kan model structure and the model structure of
Theorem~\ref{holim} is then known as the Heller model structure
\cite{MR920963}.

These two model structures are complicated to use since their
respectively cofibrant and fibrant replacement functors are not easy
to understand. The Reedy theory that is going to be exposed now is
much simpler.  This new approach allows to work with model categories
which are not necessarily cofibrantly generated: it is useful for
example for \cite{model2} Theorem~IV.3.10 where the model category
$\mathcal{M}$ is the Str{\o}m model category structure~\footnote{The
  Str{\o}m model category is conjecturally not cofibrantly generated.}
of compactly generated topological spaces with the homotopy
equivalences as weak equivalences \cite{MR35:2284} \cite{MR39:4846}
\cite{ruse} \cite{strom2}. But it requires a particular structure on
the small category $\mathcal{B}$:

\bd Let $\mathcal{B}$ be a small category. A {\rm Reedy structure} on
$\mathcal{B}$ consists of two subcategories $\mathcal{B}_{-}$ and
$\mathcal{B}_{+}$, a function $d:\mathcal{B}\longrightarrow \lambda$
called the degree function for some ordinal $\lambda$, such that every
non-identity map in $\mathcal{B}_{+}$ raises degree, every
non-identity map in $\mathcal{B}_{-}$ lowers degree, and every map
$f\in \mathcal{B}$ can be factored uniquely as $f=g\circ h$ with $h\in
\mathcal{B}_{-}$ and $g\in \mathcal{B}_{+}$. A small category together
with a Reedy structure is called a Reedy category.  \ed

Let $\mathcal{B}$ be a Reedy category. Let $b$ be an object of
$\mathcal{B}$. The \textit{latching category}
$\de(\mathcal{B}_{+}\ddownarrow b)$ at $b$ is the full
subcategory of the comma category $\mathcal{B}_{+}\ddownarrow b$
containing all the objects except the identity map of $b$. The
\textit{matching category} $\de(b\ddownarrow\mathcal{B}_{-})$ at
$b$ is the full subcategory of the comma category
$b\ddownarrow\mathcal{B}_{-}$ containing all the objects except
the identity map of $b$.

\bd Let $\mathcal{M}$ be a complete and cocomplete category. Let
$\mathcal{B}$ be a Reedy category. Let $b$ be an object of
$\mathcal{B}$. The {\rm latching space functor} is the composite
$L_b:\mathcal{M}^\mathcal{B}\longrightarrow
\mathcal{M}^{\de(\mathcal{B}_{+}\ddownarrow b)}\longrightarrow
\mathcal{M}$ where the latter functor is the colimit functor.  The
{\rm matching space functor} is the composite
$M_b:\mathcal{M}^\mathcal{B}\longrightarrow
\mathcal{M}^{\de(b\ddownarrow\mathcal{B}_{-})}\longrightarrow
\mathcal{M}$ where the latter functor is the limit functor. \ed

The \textit{Reedy model structure} of $\mathcal{M}^\mathcal{B}$ is
then constructed as follows:
\begin{itemize}
\item The \textit{Reedy weak equivalences} are the objectwise weak
  equivalences.
\item The \textit{Reedy cofibrations} are the morphisms of diagrams
  from $X$ to $Y$ such that for every object $b$ of $\mathcal{B}$ the
  morphism $X_b\sqcup_{L_b X} L_bY \longrightarrow Y_b$ is a
  cofibration of $\mathcal{M}$.
\item The \textit{Reedy fibrations} are the morphisms of diagrams from
  $X$ to $Y$ such that for every object $b$ of $\mathcal{B}$ the
  morphism $X_b\longrightarrow Y_b\p_{M_b Y} M_bX$ is a fibration of
  $\mathcal{M}$.
\end{itemize}

\bth \label{Reedymodel} (D.M. Kan) (\cite{ref_model2} Theorem~15.3.4,
Theorem~15.3.15 and Theor\-em~15.6.\-27) Let $\mathcal{B}$ be a Reedy
category. Let $\mathcal{M}$ be a model category. The objectwise weak
equivalences together with the Reedy cofibrations and the Reedy
fibrations assemble to a structure of model category. Moreover, one
has:
\begin{itemize}
\item If $\mathcal{M}$ is a cofibrantly generated for which there are
  a set $I$ of generating cofibrations whose domains and codomains are
  small relative to $I$ and a set $J$ of generating trivial
  cofibrations whose domains and codomains are small relative to $J$,
  then the Reedy model structure is cofibrantly generated.
\item If $\mathcal{M}$ is left proper (resp. right proper, proper),
  then so is $\mathcal{M}^\mathcal{B}$.
\item A morphism of diagrams from $X$ to $Y$ is a trivial Reedy
  cofibration if and only if for every object $b$ of $\mathcal{B}$ the
  morphism $X_b\sqcup_{L_b X} L_bY \longrightarrow Y_b$ is a trivial
  cofibration of $\mathcal{M}$.
\item A morphism of diagrams from $X$ to $Y$ is a trivial Reedy
  cofibration if and only if for every object $b$ of $\mathcal{B}$ the
  morphism $X_b\longrightarrow Y_b\p_{M_b Y} M_bX$ is a trivial
  fibration of $\mathcal{M}$.
\end{itemize}
\eth

Theorem~\ref{Reedymodel} is used in \cite{3eme} Theorem~8.4 and in
\cite{4eme} Theorem~8.1, Theorem~8.2 and Theorem~8.3.

The calculation of homotopy colimits is then possible in the following
situation:

\bth \label{Rcolim} (\cite{ref_model2} Proposition~15.10.2 and
Theorem~15.10.8) Let $\mathcal{B}$ be a Reedy category. Then the
following conditions are equivalent:
\begin{itemize}
\item For every fibrant object $X$ of every model category
  $\mathcal{M}$, the diagram $\cons (X)$ is Reedy fibrant (one says
  that the category $\mathcal{B}$ has {\rm fibrant constants}).
\item For every object $b$ of $\mathcal{B}$, the matching category
  $\de(b\ddownarrow\mathcal{B}_{-})$ is either empty or connected.
\item For every model category $\mathcal{M}$, the categorical
  adjunction $\liminj : \mathcal{M}^\mathcal{B} \rightleftarrows
  \mathcal{M}:\cons$ is a Quillen adjunction.
\end{itemize}
\eth

The statements of \cite{3eme} Theorem~7.5 and of \cite{4eme}
Corollary~7.4 exhibit non-trivial examples of small categories having
fibrant constants: the Reedy categories are constructed from the
category of simplices, i.e. from the order complex \cite{MR493916}, of
a locally finite poset. These constructions are reused in the proof of
\cite{ccsprecub} Theorem~7.8.

The dual statement allows the calculation of homotopy limits:

\bth \label{Rlim} (\cite{ref_model2} Proposition~15.10.2 and Theorem~15.10.8) 
Let $\mathcal{B}$ be a Reedy category. Then the following
conditions are equivalent:
\begin{itemize}
\item For every cofibrant object $X$ of every model category
  $\mathcal{M}$, the diagram $\cons (X)$ is Reedy cofibrant (one says
  that the category $\mathcal{B}$ has {\rm cofibrant constants}).
\item For every object $b$ of $\mathcal{B}$, the latching category
  $\de(\mathcal{B}_{+}\ddownarrow b)$ is either empty or connected.
\item For every model category $\mathcal{M}$, the categorical
  adjunction $\cons : \mathcal{M} \rightleftarrows
  \mathcal{M}^\mathcal{B}:\limproj$ is a Quillen adjunction.
\end{itemize}
\eth

\section{Application of the Reedy theory to the homotopy theory of flows}

The following situations are applications of Theorem~\ref{Rcolim} and
of Theorem~\ref{Rlim} to the homotopy theory of flows.

\subsection*{Homotopy pushout} 
Let $A$, $B$ and $C$ be three cofibrant objects of $\mathcal{M}$.  The
colimit of the diagram $A\longleftarrow B \longrightarrow C$ is a
homotopy colimit as soon as one of the map $B\rightarrow A$ or $B
\rightarrow C$ is a cofibration. In particular, consider an objectwise
weak equivalence $f$ from a diagram $X_1:A_1\longleftarrow B_1
\longrightarrow C_1$ to a diagram $X_2:A_2\leftarrow B_2 \rightarrow
C_2$ such that the objects $A_1, B_1, C_1, A_2, B_2, C_2$ are
cofibrant and such that the maps $B_i \rightarrow C_i$ are
cofibrations for $i=1,2$. Then $\liminj f:\liminj X_1 \rightarrow
\liminj X_2$ is a weak equivalence of $\mathcal{M}$
(\cite{MR99h:55031} Lemma~5.2.6 called the cube lemma by M. Hovey).
The proof consists of considering the Reedy category $\mathcal{I}:0
\longleftarrow 1 \longrightarrow 2$ with the degree equal to the
corresponding object. In particular, the small category $\mathcal{I}$
has fibrant constants. By Theorem~\ref{Rcolim}, the colimit functor
$\liminj : \mathcal{M}^\mathcal{I} \rightarrow \mathcal{M}$ is
therefore a left Quillen adjoint.  One has the equalities
\begin{itemize}
\item $M_0X = M_2X = \mathbf{1}$ and $M_1 X =X_0$
\item $L_0X = L_1X = \varnothing$ and $L_2X=X_1$
\item $M_0Y = M_2Y = \mathbf{1}$ and $M_1 Y =Y_0$
\item $L_0Y = L_1Y = \varnothing$ and $L_2Y=Y_1$. 
\end{itemize}
Thus, a map of diagrams $f:X\rightarrow Y$ is a Reedy cofibration if
and only if $f_0:X_0\rightarrow Y_0$ and $f_1:X_1\rightarrow Y_1$ are
cofibrations and the map $X_2\sqcup_{X_1} Y_1 \rightarrow Y_2$ is a
cofibration. So with the hypothesis above, the diagrams $X$ and $Y$
are both Reedy cofibrant. Thus, the colimit functor takes the
objectwise weak equivalence $f$ to a weak equivalence $\liminj f$ and
the colimits of the diagrams $X$ and $Y$ give their homotopy colimit.
This technique is used in the proof of \cite{exbranch} Lemma~8.5.

\subsection*{Homotopy colimit of tower of cofibrations between cofibrant objects} 
Let $(A_n\rightarrow A_{n+1})_{n\geq 0}$ be a family of cofibrations
between cofibrant objects. Then the colimit $\liminj A_n$ is a
homotopy colimit. It suffices to consider the Reedy category $0
\longrightarrow 1 \longrightarrow 2 \longrightarrow \dots$ All
matching categories are empty so this category has fibrant constants.
And $L_0A=\varnothing$, $L_nA=A_{n-1}$ for any $n\geq 1$. So the tower
is Reedy cofibrant. Hence the colimit gives the homotopy colimit.

\subsection*{Homotopy colimit of tower of cofibrations in a left proper model category} 
Let $(A_n\rightarrow A_{n+1})_{n\geq 0}$ be a family of cofibrations.
Assume $\mathcal{M}$ left proper. Then the colimit $\liminj A_n$ is a
homotopy colimit again. We do not suppose anymore the objects $A_i$
cofibrant but we must suppose $\mathcal{M}$ left proper.  The
situation is very close to the preceding situation. In fact, one is
reduced to working in the comma category $A_0\ddownarrow \mathcal{M}$
after a clever argument due to D. M. Kan using left properness (see
\cite{ref_model2} Proposition~17.9.3). This technique is used in the
proof of \cite{3eme} Theorem~11.2 and in the proof of \cite{4eme}
Theorem~9.1. What matters is the left properness of the category of
flows which is proved in \cite{2eme} Theorem~7.4.

\subsection*{Homotopy pullback} 
Let $A$, $B$ and $C$ be three fibrant objects of $\mathcal{M}$.  The
limit of the diagram $A \longrightarrow B \longleftarrow C$ is a
homotopy limit as soon as one of the maps $A\rightarrow B$ or $B
\leftarrow C$ is a fibration. This situation is used for the proof of
\cite{model2} Theorem~IV.3.14 where the model category $\mathcal{M}$
is again the Str{\o}m model category structure of compactly generated
topological spaces with the homotopy equivalences as weak equivalences
\cite{MR35:2284} \cite{MR39:4846} \cite{ruse}.

\subsection*{Homotopy limit of tower of fibrations between fibrant objects}
Let $(A_{n+1}\rightarrow A_{n})_{n\geq 0}$ be a family of fibrations
between fibrant objects. Then the limit $\limproj A_n$ is a homotopy
limit. This situation is used for the proof of \cite{model2}
Theorem~IV.3.10.

\subsection*{Homotopy limit of tower of fibrations in a right proper model category} 
Let $(A_{n+1}\rightarrow A_{n})_{n\geq 0}$ be a family of fibrations.
Assume $\mathcal{M}$ right proper. Then the limit $\limproj A_n$ is a
homotopy limit.

\subsection*{Category of simplices and of cubes and homotopy colimit}
Let $K$ be a simplicial set \cite{MR2001d:55012}. Consider the
category $\Delta K$ of simplices of $K$ defined as follows
($\Delta[n]$ being the $n$-simplex): the objects are the maps
$\Delta[n]\rightarrow K$ and the morphisms are the commutative
diagrams of simplicial sets 
\[
\xymatrix{
\Delta[m] \ar@{->}[rr]\ar@{->}[dr]&& \Delta[n] \ar@{->}[dl]\\
& K & }
\]
In other terms, $\Delta K$ is the comma category $\Delta \ddownarrow
K$ where $\Delta$ is the small category such that the presheaves over
$\Delta$ are the simplicial sets. Then

\bp \label{simp} (\cite{ref_model2} Proposition~15.10.4) Let $K$ be a
simplicial set.  The category of simplices $\Delta K$ is a Reedy
category which has fibrant constants. \ep

Similarly, let $K$ be a precubical set, that is a cubical set without
degeneracy maps \cite{Brown_cube} \cite{ccsprecub}. Let $\square$ be
the small category such that the presheaves over $\square$ are the
precubical sets. Then the category of cubes $\square \ddownarrow K$ of
the precubical set $K$ is a Reedy category which has fibrant
constants. In particular, this implies (with a little work) that the
geometric realization of a precubical set as a flow $|K| :=
\liminj_{\square[n]\rightarrow K}
(\{\widehat{0}<\widehat{1}\}^n)^{cof}$ as defined in \cite{ccsprecub}
is actually an homotopy colimit. So one has the weak S-homotopy
equivalence $|K| \simeq \holiminj _{\square[n]\rightarrow K}
\{\widehat{0}<\widehat{1}\}^n$.

\end{document}